\documentclass[12pt,a4paper]{amsart}
\usepackage{amsmath,amssymb,amscd}
\input xy
\xyoption{all}

\newcommand{\ra}{\rightarrow}

\newcommand{\PP}{\mathbb P}

\newcommand{\cO}{\mathcal{O}}

\theoremstyle{plain}
\newtheorem{theorem}{Theorem}[section]

\newtheorem{cor}[theorem]{Corollary}
\newtheorem{remark}[theorem]{Remark}

\begin{document}
\title[Clifford's theorem]{Clifford's theorem for coherent systems}

\author{H. Lange}
\author{P. E. Newstead}

\address{H. Lange\\Mathematisches Institut\\
              Universit\"at Erlangen-N\"urnberg\\
              Bismarckstra\ss e $1\frac{ 1}{2}$\\
              D-$91054$ Erlangen\\
              Germany}
              \email{lange@mi.uni-erlangen.de}
\address{P.E. Newstead\\Department of Mathematical Sciences\\
              University of Liverpool\\
              Peach Street, Liverpool L69 7ZL, UK}
\email{newstead@liv.ac.uk}
\thanks{Both authors are members of the research group VBAC (Vector Bundles on Algebraic Curves).  
         The second author would like to thank the Mathematisches Institut der Universit\"at 
         Erlangen-N\"urnberg for its hospitality}

\begin{abstract}
Let $C$ be an algebraic curve of genus $g \geq 2$. We prove an analogue of Clifford's theorem for coherent
systems on $C$ and some refinements using results of Re and Mercat.
\end{abstract}
\maketitle

\section{Introduction}

Let $C$ be a smooth projective algebraic curve of genus $g \geq 2$ defined over an algebraically closed field. The 
classical Clifford's theorem states that for any line bundle $L$ on $C$ of degree $d$, $0 \leq d \leq 2g-2$, we have
$$
h^0(L) \leq \frac{d}{2} + 1.
$$
Moreover one can define the Clifford index $\gamma_C$ of $C$ by
$$
\gamma_C := \mbox{Min} \{ d - 2h^0(L) + 2\; |\; \deg L = d,\; h^0(L) \geq 2,\; h^1(L) \ge 2 \}.
$$
This definition works for $g \geq 4$, since there always exist line bundles with this property. For $g=2$ and for 
hyperelliptic curves of genus 3 we define $\gamma_C = 0$ and for non-hyperelliptic curves of genus 3 we define 
$\gamma_C = 1$. For low values of $\gamma_C$ a complete description of curves of Clifford index $\gamma_C$ is known.
In particular, $\gamma_C = 0$ if and only if $C$ is hyperelliptic. Note also that for line bundles outside the range of 
Clifford's theorem either $h^0(L) \leq 1$ or $h^1(L) \leq 1$ and the Riemann-Roch theorem then gives a 
bound for $h^0(L)$.\\

Clifford's theorem has been generalized to semistable vector bundles. In fact, for any semistable vector bundle $E$ of 
rank $n$ and degree $d$ with $0 \leq d \leq (2g-2)n$,
\begin{equation} \label{e1}
h^0(E) \leq \frac{d}{2} + n
\end{equation}
(see \cite[Theorem 2.1]{Bgn}). Moreover refinements of this for curves of Clifford index $\geq 1$ and for Clifford 
index $\geq 2$ have been obtained by Re \cite{re} and Mercat \cite{me}. It is again easy to deduce bounds for 
all values of $d$, not just those in the range to which Clifford's theorem and its refinements apply. However
no good generalization to $\alpha$-stable coherent systems is known. Of course, if $(E,V)$ is $\alpha$-stable 
for small positive values of $\alpha$ (or $0^+$-stable in the notation of \cite{bgvn}), 
then $E$ is a semistable  bundle and \eqref{e1} applies, but there is no a priori reason why this should 
extend to higher values of $\alpha$.

In fact, there are two questions which can be asked here. The first question asks for a good bound on $h^0(E)$  
for any $\alpha$-semistable coherent system. The fact that this bound exists is known; indeed a bound 
(possibly dependent on $\alpha)$ is needed in order to construct the moduli spaces of semistable coherent systems
as quasi-projective varieties (see for example \cite[Corollary 2.6.1]{kn} 
and for a more precise bound \cite[Lemmas 10.2 and 10.3]{bgvn}). Once one has such a bound, it follows from
\cite[Propositions 4.2 and 4.6]{bgvn} that, for given $(n,d)$, one needs to test the bounds for only a finite number of values 
of $\alpha$. So there exists a  bound independent of $\alpha$. However one expects a much better bound than those that are known. 

The second question is to determine for given $(n,d)$ the maximum value of $k$ for which there exists an $\alpha$-stable 
or $\alpha$-semistable coherent system $(E,V)$ of type $(n,d,k)$. Note that $k$ could be $< h^0(E)$, so the two 
questions are different.\\

Our main object in this paper is to address the second question. We shall prove in fact the following analogue 
of Clifford's theorem, where for completeness  we include a result of Riemann-Roch type. \\

\noindent
{\bf Theorem 2.1}. 
{\it Let $(E,V)$ be a coherent system of type $(n,d,k)$ with $k > 0$ which is $\alpha$-semistable for some 
$\alpha > 0$ on a smooth 
projective curve $C$ of genus $g \geq 2$. Then
\begin{eqnarray*}  
k &\leq &\frac{d}{2} + n   \qquad \mbox{if} \qquad 0 \leq d \leq 2gn  \quad \mbox{and}\\
k &\leq & d + n(1-g)        \quad \mbox{if} \qquad d \geq 2gn.
\end{eqnarray*}}

The proof depends not only on the result for semistable bundles, but also on some parts of the theory 
developed in \cite{bgvn}. As consequences of a general theorem (Theorem \ref{thm3.1}) which may have further 
applications, we obtain refinements for the cases $\gamma_C \geq 1$ and $\gamma_C \geq 2$
(see Corollaries \ref{cor3.2} and \ref{cor3.3}). We should remark that the Clifford index as defined above 
will not necessarily give precise bounds; a refinement of the definition of the Clifford index is certainly 
needed for coherent systems of rank $n \geq 2$.

We include a complete list of cases for which the Clifford inequality in Theorem \ref{thm2.1} is an equality 
(Remark \ref{rem2.2}) and an example to show that $h^0(E)$ 
does not necessarily satisfy either inequality (Remark \ref{rem2.3}).

For curves of genus 0 and 1 Clifford's theorem makes no sense. However in these cases general bounds on $k$ can be obtained 
from the results of \cite{ln1} and \cite{ln2} (see Remark \ref{rem2.4}).

\section{Clifford's theorem}

\begin{theorem} \label{thm2.1}
Let $(E,V)$ be a coherent system of type $(n,d,k)$ with $k > 0$ which is $\alpha$-semistable for some 
$\alpha > 0$ on a smooth 
projective curve $C$ of genus $g \geq 2$. Then
\begin{eqnarray}  \label{eqn1}
k &\leq &\frac{d}{2} + n   \qquad \mbox{if} \qquad 0 \leq d \leq 2gn
\end{eqnarray}
and
\begin{eqnarray}  \label{e3}
k &\leq & d + n(1-g)        \quad \mbox{if} \qquad d \geq 2gn.
\end{eqnarray}
\end{theorem}

\begin{proof}
The proof is by induction on $n$, the case $n=1$ being the classical theorems of Clifford and Riemann-Roch. So suppose $n \geq 2$ and 
the theorem is proved for coherent systems of rank less than $n$.

If $(E,V)$ is $0^+$-stable, then the bundle $E$ itself is semistable and hence the assertion holds by the known 
Clifford's theorem for semistable bundles (see \cite[Theorem 2.1]{Bgn}) and the Riemann-Roch theorem. (Note 
that for a semistable bundle $E$ of degree $> (2g-2)n$ we have $h^1(E) = 0$.)

If $(E,V)$ is not $0^+$-stable, then it must be strictly $\alpha$-semistable for some $\alpha >0$ (see \cite[Lemma 6.5]{bgvn}).
So for this $\alpha$, there exists 
an exact sequence
\begin{equation}  \label{eqn2}
0 \ra (E_1,V_1) \ra (E,V) \ra (E_2,V_2) \ra 0
\end{equation}
with $(E_i,V_i)$ of type $(n_i,d_i,k_i)$ for $i = 1, 2$ and
\begin{equation} \label{eq}
\frac{d_1}{n_1} + \alpha \frac{k_1}{n_1} = \frac{d_2}{n_2} + \alpha \frac{k_2}{n_2}.
\end{equation} 
Hence $(E_1,V_1)$ and $(E_2,V_2)$ are $\alpha$-semistable. Moreover, we can assume that $(E_1,V_1)$ contradicts 
the $0^+$-stability of $(E,V)$ and so
\begin{equation} \label{eq4} 
\frac{d_1}{n_1} \geq \frac{d_2}{n_2} \quad \mbox{and} \quad \frac{k_1}{n_1} \leq \frac{k_2}{n_2}.
\end{equation} 
It follows that $k_2 > 0$ and so $d_2 \geq 0$. 

Suppose first $0 \leq d \leq 2gn$. Then clearly $\frac{d_2}{n_2} \leq 2g$.
By the induction hypothesis we get
\begin{equation} \label{e6}
k_2 \leq \frac{d_2}{2} + n_2.
\end{equation}
On the other hand,
$$
\frac{k_1}{n_1} \leq \frac{k_2}{n_2} \leq \frac{d_2}{2n_2} + 1 \leq \frac{d_1}{2n_1} + 1
$$ 
implying
\begin{equation}  \label{e7}
k_1 \leq \frac{d_1}{2} + n_1.
\end{equation}
Adding \eqref{e6} and \eqref{e7} we get \eqref{eqn1}.

Now suppose $d \geq 2gn$. If $k \leq \frac{d}{2} + n$, then the result holds. 
Otherwise $k > \frac{d}{2} + n$. Hence $k_2 > \frac{d_2}{2} + n_2$ by \eqref{eq4} and $d_2 > 2gn_2$ by the induction hypothesis
and \eqref{eqn1}. So again by the induction hypothesis,
\begin{equation} \label{e8}
k_2 \leq d_2 + n_2(1-g).
\end{equation}
It follows that
$$
\frac{k_1}{n_1} \leq \frac{k_2}{n_2} \leq \frac{d_2}{n_2} + 1 -g \leq \frac{d_1}{n_1} + 1 - g,
$$
i.e.  
\begin{equation} \label{e9}
k_1 \leq d_1 + n_1(1-g).
\end{equation}
Adding \eqref{e8} and \eqref{e9} gives \eqref{e3}. 
\end{proof}

\begin{remark} \label{rem2.1}
{\em If $d < 0$, there do not exist $\alpha$-semistable coherent systems for any $k > 0$ 
(see \cite[Lemmas 4.1 and 4.3]{bgvn})}.
\end{remark}

\begin{remark} \label{rem2.2}
{\em If $(E,V)$ is an $\alpha$-semistable coherent system for which the underlying bundle $E$ is not semistable, 
then we can choose the sequence \eqref{eqn2} so that the inequalities in \eqref{eq4} are strict. Hence so 
are the inequalities of the theorem. On the other hand, if $d \geq 2gn$ and $(E,V) = (E,H^0(E))$ with $E$
semistable, then equality holds in \eqref{e3} by the Riemann-Roch theorem. 

If $C$ is non-hyperelliptic and $d < 2gn$, 
then equality in \eqref{eqn1} holds only for $(E,V) \simeq (\cO^{\oplus n}, H^0(\cO^{\oplus n}))$ and
$(E,V) \simeq (K^{\oplus n}, H^0(K^{\oplus n}))$. These coherent systems are strictly $\alpha$-semistable for 
all $\alpha >0$ if $n \geq 2$. If $C$ is hyperelliptic, we must add to this list 
$(E,V) \simeq ((L^a)^{\oplus n}, H^0((L^a)^{\oplus n}))$, where $L$ is the hyperelliptic line bundle and $1 \leq a \leq g-2$.
Again these are strictly $\alpha$-semistable for all $\alpha > 0$ if $n \geq 2$. Thus, for any $\alpha$-stable coherent 
system with $n \geq 2$ and $0 \leq d < 2gn$ the inequality \eqref{eqn1} is strict.
}
\end{remark}

\begin{remark} \label{rem2.3}
{\em It is easy to give examples of $\alpha$-stable coherent systems $(E,V)$ such that $h^0(E)$ does not satisfy 
the bounds of Theorem \ref{thm2.1}.

To give an example, consider an exact sequence
$$
0 \ra (L_1,V_1) \ra (E,V) \ra (L_2,V_2) \ra 0,
$$ 
where $L_1$ and $L_2$ are line bundles of degrees $d_1$ and $d_2$ with $d_1 > d_2$ and $\dim V_1 = 1, \; \dim V_2 = 2$.
Then
$$
\dim \mbox{Ext}^1((L_2,V_2),(L_1,V_1)) \geq d-g-1
$$
by \cite[equations (8) and (9)]{bgvn}. So we can choose the extension to be non-trivial provided $d \geq g + 2$. Then, by the 
results of \cite[section 6]{bgvn}, the coherent system $(E,V)$ is $\alpha_c^+$-stable for $\alpha_c = d_1 - d_2$.
On the other hand,
\begin{equation} \label{eq5}
h^0(E) \geq 2 + h^0(L_1) \geq d_1 - g + 3
\end{equation}
by Riemann-Roch. Now suppose for simplicity that $g$ is even. Then by classical Brill-Noether theory we can take 
$d_2 = \frac{g}{2} + 1$. So \eqref{eq5} says $$
h^0(E) \geq d - \frac{3g}{2} + 2,
$$ 
which is larger than the estimate for $k$ of Theorem \ref{thm2.1}, provided that $d > 3g$.}
\end{remark}

\begin{remark}  \label{rem2.4}
{\em The above theorem does not apply for $g=0$ and $g=1$ and of course Clifford's theorem does not exist in these cases. 
However, we do know that for any semistable coherent system $(E,V)$ on $\PP^1, \; \dim V \leq  d+n$ 
(see \cite[Lemma 3.1]{ln1}) and in fact $h^0(E) = d+n$. 

On a curve of genus 1, if $(E,V)$ is $\alpha$-stable 
and $n \geq 2$, then $\dim V \leq d$ (see \cite[Corollary 4.2]{ln2}). 
If $E$ is semistable of degree $d > 0$, then $E$ is a direct sum of indecomposable bundles, all of positive degree, 
so $h^0(E) = d$ by \cite[Lemma 15]{at}. If $E$ is semistable of degree 0, then $h^0(E) \leq n$ with equality if and 
only if $E \simeq \cO^{\oplus n}$; moreover $(\cO^{\oplus n}, H^0(\cO^{\oplus n}))$ is $\alpha$-semistable for all 
$\alpha > 0$. Finally, if $(E,V)$ is $\alpha$-semistable for some $\alpha > 0$ and $E$ is not semistable, then $(E,V)$ 
must contain a subsystem $(\cO, H^0(\cO))$ (see proof of \cite[Lemma 4.1]{ln1}) and one can deduce that 
$\dim V \leq \max \{n-1,d \}$. The equality $\dim V = n-1$ is attained for $1 \leq d < n-1$ by
$$
(E,V) \simeq (\cO^{\oplus r}, H^0(\cO^{\oplus r})) \oplus (G,W),
$$
where $1 \leq r \leq n-1, \; \deg G = d$ and $\dim W = n-r-1$. This coherent system is $\alpha$-semistable for 
$\alpha = d$ provided that $(G,W)$ is $\alpha$-semistable for the same value of $\alpha$. By \cite[Theorem 5.2]{ln1},
such $\alpha$-semistable coherent systems $(G,W)$ exist whenever $d \geq n-r-1$.  
}
\end{remark}

\section{Refinements of Clifford's theorem}

Our main purpose in this section is to show that the refinements of Clifford's theorem for semistable bundles 
in \cite{re} and \cite{me} can be extended to $\alpha$-semistable coherent systems. In fact, we prove something 
more general which could have applications for curves with Clifford index $\gamma_C \geq 3$.\\ 

Let $\delta_1, \ldots, \delta_r, a_1, \ldots , a_{r-1}, b_1, \ldots, b_{r-1}$ be rational numbers such that
\begin{itemize}
\item $1 = \delta_1 < \delta_2 < \cdots < \delta_r$,
\item $0 < a_1 < \cdots < a_{r-1}$,
\item $a_1 + b_1 = 1$,
\item$a_i \delta_{i+1} + b_i = a_{i+1} \delta_{i+1} + b_{i+1}$  for $1 \leq i \leq r-2$.
\end{itemize}
Note that the conditions imply that $a_i \mu + b_i > 1$ whenever $\delta_i < \mu\le\delta_{i+1}$.

Consider the following condition on $k$ and $n$:
\begin{equation} \label{eq6}
\delta_i \leq \mu \leq \delta_{i+1}  \Rightarrow k \leq (a_i\mu + b_i) n \quad \mbox{for} \quad 1 \leq i \leq r-1.
\end{equation}

\begin{theorem} \label{thm3.1}
Suppose $k \geq 1$ and that condition \eqref{eq6} holds for any semistable bundle $E$ 
of rank $n$ and degree $d$ with $\frac{d}{n} = \mu$ and $h^0(E) = k$. 
Then the same condition holds for any coherent
system $(E,V)$ which is $\alpha$-semistable for some $\alpha > 0$ and of type $(n,d,k)$ with $\frac{d}{n} = \mu$. 
\end{theorem}

\begin{proof}
The proof is by induction on $n$, the case $n=1$ being immediate. So suppose $n \geq 2$ and the theorem is proved for 
coherent systems of rank less than $n$.

Let $(E,V)$ be of type $(n,d,k)$. If $(E,V)$ is $0^+$-stable, then $E$ is semistable and the result is clear.
Otherwise, as in the proof of Theorem \ref{thm2.1}, we have an exact sequence \eqref{eqn2} with conditions \eqref{eq} 
and \eqref{eq4}.
Suppose first that $\delta_1 \leq \mu \leq \delta_2$. If $k \leq n$, the condition \eqref{eq6} holds for any $\mu \geq 1$.
So suppose $k > n$. Then $k_2 > n_2$ by \eqref{eq4} and $\frac{d_2}{n_2} > 1$ by \cite[Lemma 3.1]{bgmn}; moreover
$\frac{d_2}{n_2} \leq  \frac{d}{n} \leq \delta_2$ by \eqref{eq4}. So by the induction hypothesis,
\begin{equation} \label{ee1}
k_2 \leq a_1d_2 + b_1n_2
\end{equation}
and hence also
$$
\frac{k_1}{n_1} \leq \frac{k_2}{n_2} \leq a_1 \frac{d_2}{n_2} + b_1 \leq a_1 \frac{d_1}{n_1} + b_1,
$$
i.e.
\begin{equation} \label{ee2}
k_1 \leq a_1d_1 + b_1n_1
\end{equation}
Adding \eqref{ee1} and \eqref{ee2} gives $k \leq a_1d + b_1n$ as required. 

Now suppose $r \geq 3$ and the result holds for $1 \leq \mu \leq \delta_{s-1}$ with $3 \leq s \leq r$. 
Suppose $\delta_{s-1} \leq \mu \leq \delta_s$. 
Then,
if $k \leq (a_{s-1}\delta_{s-1} + b_{s-1})n$, the result holds. So suppose $k > (a_{s-1}\delta_{s-1} + b_{s-1})n$.
Then $k_2 > (a_{s-1}\delta_{s-1} + b_{s-1})n_2$ and so $\frac{d_2}{n_2} > \delta_{s-1}$.
Clearly $\frac{d_2}{n_2} \leq \frac{d}{n} \leq \delta_s$. Hence by the induction hypothesis
\begin{equation} \label{ee3}
k_2 \leq a_{s-1}d_2 + b_{s-1} n_2.
\end{equation}
Since $\frac{k_1}{n_1} \leq \frac{k_2}{n_2}$ and $\frac{d_1}{n_1} \geq \frac{d_2}{n_2}$ we get
\begin{equation} \label{ee4}
k_1 \leq a_{s-1} d_1 + b_{s-1} n_1.
\end{equation}
Adding \eqref{ee3} and \eqref{ee4} gives $k \leq a_{s-1} d + b_{s-1}n$ for $\delta_{s-1} \leq \mu \leq \delta_s$. 
This completes the proof of the theorem by induction.
\end{proof}

\begin{remark} \label{rem3.2}
{\em For $0 \leq \frac{d}{n} < 1$ the best inequality of the above type we can give is $k \leq n$.
However, if we assume $0 < \frac{d}{n} < 1$, then
$$
k \leq n + \frac{1}{g}(d-n)
$$
for any $\alpha$-semistable coherent system with $\alpha > 0$ (see \cite[Lemma 3.1]{bgmn}). Moreover this 
inequality is precisely the condition for the existence of $\alpha$-stable coherent systems in this range.
}
\end{remark}

We can now extend Re's result (see \cite[Propositions 3 and 4]{re}) to coherent systems.

\begin{cor} \label{cor3.2}
Let $(E,V)$ be a coherent system of type $(n,d,k)$ with $k > 0$ which is $\alpha$-semistable for some $\alpha > 0$ 
on a non-hyperelliptic curve $C$ of genus $g \geq 3$.  \\
{\em (a):} Suppose $n \leq d \leq (2g-3)n$. Then
$$
k \leq \frac{d+n}{2}.
$$
{\em (b):} Suppose $(2g-3)n \leq d \leq (2g-2)n$. Then
$$
k \leq d + (2-g)n.
$$
\end{cor}

\begin{proof}
Take $r=3, \; \delta_1 = 1, \; \delta_2 = 2g-3, \; \delta_3 = 2g-2$ and $a_1 = b_1 = \frac{1}{2}, a_2 = 1, b_2 = 2-g$.
Then the hypothesis of Theorem \ref{thm3.1}  follows from \cite[Propositions 3 and 4]{re}. Now the assertion is a special case of 
the theorem. 
\end{proof}

The proof implies that we have strict inequality in (a) unless $E$ is a semistable bundle.

We have also an extension of Mercat's result (see \cite[Theorem 2.1]{me}) for curves of Clifford index $\gamma_C \geq 2$.

\begin{cor} \label{cor3.3}
Let $(E,V)$ be a coherent system of type $(n,d,k)$ with $k > 0$ which is $\alpha$-semistable for some $\alpha > 0$ 
on a curve $C$ of genus $g$ with Clifford index $\gamma_C \geq 2$.\\
{\em (a):} Suppose $2 + \frac{2}{g-4} \leq \frac{d}{n} \leq 2g-4- \frac{2}{g-4}$. Then
$$
k \leq \frac{d}{2}.
$$
{\em (b):} Suppose $1 \leq \frac{d}{n} \leq 2 + \frac{2}{g-4}$. Then
$$
k \leq \frac{1}{g-2}(d-n) +n.
$$
\end{cor}

\begin{proof}
Note first that $\gamma_C \geq 2$ implies $g \geq 5$. Take $r=3, \; \delta_1 = 1, \; \delta_2 = 2 + \frac{2}{g-4}, \;
\delta_3 = 2g-4 - \frac{2}{g-4}$ and $ a_1 = \frac{1}{g-2}, \; b_1 = \frac{g-3}{g-2}, \; a_2 = \frac{1}{2}, \; b_2 = 0$.
Then
$$
\frac{1}{g-2}\left(2 + \frac{2}{g-4}\right) + \frac{g-3}{g-2} = \frac{1}{2}\left(2 + \frac{2}{g-4}\right).
$$
Now the hypotheses of Theorem \ref{thm3.1} are given by \cite[Theorem 2.1]{me}. Hence Corollary \ref{cor3.3} is a special case of the
theorem.
\end{proof}

\end{document}